\documentclass[12pt,fleqn]{iopart}
\usepackage[english]{babel}
\usepackage{amssymb}
\usepackage{amsfonts}
\usepackage{euscript}

\begin{document}

\title[Zero-Curvature Representations for Second Heavenly Equation]{%
Contact Integrable Extensions and Zero-Curvature Representations for  the
Second Heavenly Equation
}

\author{Oleg I. Morozov}

\address{Institute of Mathematics and Statistics, University of Troms\o, Troms\o 
\, 9037, Norway
\\
Oleg.Morozov{\symbol{64}}uit.no}

\begin{abstract}
The method of contact integrable extensions is used to find new zero-curvature representation for Pleba\~nski's second heavenly equation.
\end{abstract}

\ams{58H05, 58J70, 35A30}

\maketitle

\section{Introduction}
The second Pleba\~nski's heavenly equation, \cite{Plebanski},
\begin{equation}
u_{xz}=u_{ty}+u_{yy}\,u_{zz}-u_{yz}^2,
\label{PlebanskiEquation}
\end{equation}
describes self-dual gravitational fields. This equation can be obtained  as the com\-pa\-ti\-bi\-li\-ty condition for the following system of {\sc pde}s, \cite{Husain,BogdanovKonopelchenko}, cf. \cite[Eq. (3.13)]{Plebanski},
\begin{equation}
\left\{
\begin{array}{l}
v_t = (u_{yz} + \lambda) \, v_z - u_{zz}\,v_y,
\\
v_x = u_{yy} \, v_z - (u_{yz} - \lambda)\,v_y
\end{array}
\right.
\label{covering}
\end{equation}
with an arbitrary constant $\lambda$. This condition is equivalent to the commutativity of four infinite-dimensional vector fields
\begin{eqnarray}
\widetilde{D}_t &=& \bar{D}_t 
+ \sum \limits_{i,j \ge 0} \widetilde{D}^i_y \widetilde{D}^j_z \left(
(u_{yz} + \lambda) \, v_{0,1} - u_{zz}\,v_{1,0}
\right)\,\frac{\partial}{\partial v_{i,j}},
\nonumber
\\
\widetilde{D}_x &=& \bar{D}_x 
+ \sum \limits_{i,j \ge 0} \widetilde{D}^i_y \widetilde{D}^j_z \left(
u_{yy} \, v_{0,1} - (u_{yz} - \lambda)\,v_{1,0}
\right)\,\frac{\partial}{\partial v_{i,j}},
\nonumber
\\
\widetilde{D}_y &=& \bar{D}_y 
+ \sum \limits_{i,j \ge 0} v_{i+1,j}\,\frac{\partial}{\partial v_{i,j}},
\nonumber
\\
\widetilde{D}_z &=& \bar{D}_z 
+ \sum \limits_{i,j \ge 0} v_{i,j+1}\,\frac{\partial}{\partial v_{i,j}},
\nonumber
\end{eqnarray}
where $\bar{D}_t$, $\bar{D}_x$, $\bar{D}_y$ and $\bar{D}_z$ are restrictions of the total derivatives $D_t$, $D_x$, $D_y$ and $D_z$ to the infinite prolongation of Eq. (\ref{PlebanskiEquation}).  This construction is called a {\it differential covering}, \cite{KV84} -- \cite{KV99}, or zero-curvature representation. Dually Eqs. (\ref{covering}) can be defined by means of differential 1-form 
\begin{equation}
\fl
\omega =
dv
+(v_{zz}\,v_{y}- (u_{yz}+\lambda)\,v_z)\,dt
+((u_{yz}-\lambda)\,v_{y}-u_{yy}\,v_{z})\,dx
-v_y\,dy - v_z\,dz
\label{WEform}
\end{equation}
called the {\it Wahlquist--Estabrook form} of the covering, \cite{DoddMorris}.
In \cite{Morozov2009b} we show that this form can be inferred from  a linear combination of Maurer--Cartan forms of the contact sym\-met\-ry pseudo-group of Eq. (\ref{PlebanskiEquation}). In this paper we apply to (\ref{PlebanskiEquation}) the technique of contact integrable extensions ({\sc cie}s) proposed in \cite{Morozov2009a}. We find {\sc cie}s of the structure equations of the contact symmetry pseudo-group of Eq. (\ref{PlebanskiEquation}). The analysis of these {\sc cie}s splits into two cases. In the first case integration of the {\sc cie} gives Eqs. (\ref{covering}), while in the second case we obtain new covering of the second heavenly equation.

\section{Symmetry pseudo-group of the second heavenly equation}

Let $\pi \colon \mathbb{R}^n \times \mathbb{R} \rightarrow \mathbb{R}^n$ be a vector bundle with the local base coordinates $(x^1,...,x^n)$ and the local fibre coordinate $u$; then denote by $J^2(\pi)$  the bundle of the second-order jets of sections of $\pi$, with the local coordinates $(x^i,u,u_i,u_{ij})$, $i,j\in\{1,...,n\}$, $i \le j$. For every local section $(x^i,f(x))$ of $\pi$, denote by $j_2(f)$ the corresponding 2-jet $(x^i,f(x),\partial f(x)/\partial x^i,\partial^2 f(x)/\partial x^i\partial x^j)$. A differential 1-form $\vartheta$ on $J^2(\pi)$ is called a {\it contact form} if it is annihilated by all 2-jets of local sections: $j_2(f)^{*}\vartheta = 0$. In the local coordinates every contact 1-form is a linear combination of the forms  $\vartheta_0 = du - u_{i}\,dx^i$,
$\vartheta_i = du_i - u_{ij}\,dx^j$, $i, j \in \{1,...,n\}$, $u_{ji} = u_{ij}$ (here and later we assume the summation convention, so $u_i\,dx^i = \sum_{i=1}^{n}\,u_i\,dx^i$, etc.) A local diffeomorphism $\Delta \colon J^2(\pi) \rightarrow J^2(\pi)$, 
$\Delta \colon (x^i,u,u_i,u_{ij}) \mapsto (\check{x}^i,\check{u},\check{u}_i,\check{u}_{ij})$, is called a {\it contact transformation} if for every contact 1-form $\check{\vartheta}$ the form $\Delta^{*}\check{\vartheta}$ is also contact. We denote by $\mathrm{Cont}(J^2(\pi))$  the pseudo-group  of contact transformations on $J^2(\pi)$.

Let $\EuScript{H} \subset \mathbb{R}^{(2 n+1)(n+3)(n+1)/3}$ be an open set with local coordinates $a$, $b^i_k$, $c^i$, $f^{ik}$, $g_i$, $s_{ij}$, $w^k_{ij}$, $u_{ijk}$, $i,j,k \in \{1,...,n\}$,  such that $a\not =0$, $\det (b^i_k) \not = 0$, $f^{ik}=f^{ki}$, 
$u_{ijk}=u_{jik}= u_{ikj}$. Let $(B^i_k)$ be the inverse matrix for the matrix $(b^k_l)$, so $B^i_k\,b^k_l = \delta^i_l$. We consider the {\it lifted coframe}
\[
\Theta_0 = a\, \vartheta_0,
\quad
\Theta_i = g_i\,\Theta_0 + a\,B_i^k\,\vartheta_k,
\quad
\Xi^i =c^i\,\Theta_0+f^{ik}\,\Theta_k+b_k^i\,dx^k,
\]
\begin{equation}
\Theta_{ij} = a\,B^k_i\, B^l_j\,(du_{kl}-u_{klm}\,dx^m)+s_{ij}\,\Theta_0+w_{ij}^{k}\,\Theta_k,
\label{LCF}
\end{equation}
\noindent
$i \le j$, defined on $J^2(\pi)\times\EuScript{H}$. As it is shown in \cite{Morozov2006}, the forms (\ref{LCF}) are Maurer--Cartan forms for $\mathrm{Cont}(J^2(\pi))$, that is, a local diffeomorphism
$\widehat{\Delta} \colon J^2(\pi) \times \EuScript{H} \rightarrow J^2(\pi) \times \EuScript{H}$
satisfies the conditions
$\widehat{\Delta}^{*}\, \check{\Theta}_0 = \Theta_0$,
$\widehat{\Delta}^{*}\, \check{\Theta}_i = \Theta_i$,
$\widehat{\Delta}^{*}\, \check{\Xi}^i = \Xi^i$,
and $\widehat{\Delta}^{*}\, \check{\Theta}_{ij} = \Theta_{ij}$
whenever it is projectable on $J^2(\pi)$, and its projection
$\Delta \colon J^2(\pi) \rightarrow J^2(\pi)$ is a contact transformation.

The structure equations for $\mathrm{Cont}(J^2(\pi))$ read
\begin{eqnarray}
d \Theta_0 &=& \Phi^0_0 \wedge \Theta_0 + \Xi^i \wedge \Theta_i,
\nonumber
\\
d \Theta_i &=& \Phi^0_i \wedge \Theta_0 + \Phi^k_i \wedge \Theta_k
+ \Xi^k \wedge \Theta_{ik},
\nonumber
\\
d \Xi^i &=& \Phi^0_0 \wedge \Xi^i -\Phi^i_k \wedge \Xi^k
+\Psi^{i0} \wedge \Theta_0
+\Psi^{ik} \wedge \Theta_k,
\nonumber
\\
d \Theta_{ij} &=& \Phi^k_i \wedge \Theta_{kj} - \Phi^0_0 \wedge \Theta_{ij}
+ \Upsilon^0_{ij} \wedge \Theta_0
+ \Upsilon^k_{ij} \wedge \Theta_k +   \Xi^k \wedge \Theta_{ijk},
\nonumber
\end{eqnarray}
\noindent
where the additional forms $\Phi^0_0$, $\Phi^0_i$, $\Phi^k_i$, $\Psi^{i0}$, $\Psi^{ij}$,
$\Upsilon^0_{ij}$, $\Upsilon^k_{ij}$, and $\Theta_{ijk}$ depend on differentials of the coordinates of $\EuScript{H}$.

Suppose $\EuScript{E}$ is a second-order differential equation in one dependent and $n$ independent variables. We consider $\EuScript{E}$ as a submanifold in $J^2(\pi)$. Let $\mathrm{Cont}(\EuScript{E})$ be the group of contact symmetries for $\EuScript{E}$. It consists of all the contact transformations on $J^2(\pi)$ mapping $\EuScript{E}$ to itself. 
Let $\iota_0 \colon \EuScript{E} \rightarrow J^2(\pi)$ be an embedding and 
$\iota = \iota_0 \times \mathrm{id} \colon 
\EuScript{E}\times \EuScript{H} \rightarrow J^2(\pi)\times \EuScript{H}$. Maurer--Cartan forms of the pseudo-group $\mathrm{Cont}(\EuScript{E})$ can be obtained from the forms  $\theta_0 = \iota^{*} \Theta_0$, $\theta_i= \iota^{*}\Theta_i$, $\xi^i = \iota^{*}\Xi^i$ and $\theta_{ij}=\iota^{*}\Theta_{ij}$
by means of \'Elie Cartan's method of equivalence, \cite{Cartan1}--\cite{Cartan4}, \cite{Gardner,Kamran,Olver95}, see details and examples in \cite{FelsOlver}, \cite{Morozov2002}---\cite{Morozov2009b}.

Using this method,  we find the Maurer--Cartan forms and their structure equations for the symmetry pseudo-group of Eq. (\ref{PlebanskiEquation}). The structure equations  have the following form:
\[\fl
d \theta _{0}=
\eta _{5} \wedge\theta _{0}
+ \xi^{1} \wedge\theta _{1}
+ \xi^{2} \wedge\theta _{2}
+ \xi^{3} \wedge\theta _{3} 
+ \xi^{4} \wedge\theta _{4},
\]
\[\fl
d \theta _{1}=
(\eta _{5} - \eta _{1}) \wedge\theta _{1} 
 - \eta _{3} \wedge\theta _{2}
- \theta _{44} \wedge\theta _{3} 
 + \eta _{6} \wedge\theta _{4}
 + \xi^{1} \wedge\theta _{11} 
 + \xi^{2} \wedge\theta _{12} 
+ \xi^{3} \wedge\theta _{13}
\] \[\fl \quad \quad \quad
 + \xi^{4} \wedge\theta _{14}, 
\]
\[\fl
d \theta _{2}=
- \eta _{2} \wedge\theta _{1} 
+ (\eta _{5}-\eta_4) \wedge\theta _{2} 
 + (\eta _{6} - 2\,\theta _{34})\wedge\theta _{3} 
+ \theta _{33} \wedge\theta _{4}
+ \xi^{1} \wedge\theta _{12} 
+ \xi^{2} \wedge\theta _{22}
\] \[\fl \quad \quad \quad
+ \xi^{3} \wedge\theta _{23} 
+\xi^{4} \wedge\theta _{13}, 
\]
\[\fl
d \theta _{3}=
{\textstyle \frac {1}{3}} \,(\eta _{1} -2\,\eta_4-2\,\eta_5) \wedge\theta _{3}
 - \eta _{2} \wedge\theta _{4} 
+ \xi^{1} \wedge\theta _{13} 
+\xi^{2} \wedge\theta _{23} 
+ \xi^{3} \wedge\theta _{33} 
+ \xi^{4} \wedge\theta _{34} ,
\]
\[\fl
d \theta _{4}=
- \eta _{3} \wedge\theta _{3} 
+ {\textstyle \frac {1}{3}} \,(\eta _{4} -2\,\eta_1+2\,\eta_5)\wedge\theta _{4}
+ \xi^{1} \wedge\theta _{14} 
+\xi^{2} \wedge\theta _{13} 
+ \xi^{3} \wedge\theta _{34} 
+ \xi^{4} \wedge\theta _{44}, 
\]
\[\fl
d \xi^{1}=\eta _{1} \wedge\xi^{1} + \eta _{2} \wedge\xi^{2},
\]
\[\fl
d \xi^{2}=\eta _{3} \wedge\xi^{1} + \eta _{4} \wedge\xi^{2},
\]
\[\fl
d \xi^{3}=
\theta _{44} \wedge\xi^{1} 
- (\eta _{6} - 2\,\theta _{34}) \wedge\xi^{2} 
+{\textstyle \frac {1}{3}} \,(\eta _{5} - \eta _{1} + 2 \,\eta _{4}) \wedge\xi^{3} 
+ \eta _{3} \wedge\xi^{4}, 
\]
\[\fl
d \xi^{4}= 
- \eta _{6} \wedge\xi^{1} 
- \theta _{33} \wedge\xi^{2}
+ \eta _{2} \wedge\xi^{3} 
+ {\textstyle \frac {1}{3}} \,(\eta _{5} +2\,\eta_1- \eta _{4} ) \wedge\xi^{4}, 
\]
\[\fl
d \theta _{11}= 
- \eta _{12} \wedge\theta _{1} 
+ \eta _{14} \wedge\theta _{2} 
+ \eta _{13} \wedge\theta _{3}
- \eta _{11} \wedge\theta _{4} 
+ \eta _{5} \wedge\theta _{11} 
- 2\,\eta _{1} \wedge\theta _{11} 
- 2\,\eta _{3} \wedge\theta _{12} 
\] \[\fl \quad \quad \quad
- 2\,\theta _{44} \wedge\theta _{13} 
+ 2\,\eta _{6} \wedge\theta _{14} 
- \xi^{1} \wedge\eta _{22} 
- \xi^{2} \wedge\eta _{21} 
- \xi^{3} \wedge\eta _{17} 
- \xi^{4} \wedge\eta _{18}, 
\]
\[\fl
d \theta _{12}= 
- \eta _{10} \wedge\theta _{1} 
+ \eta _{12} \wedge\theta _{2} 
+ \eta _{11} \wedge\theta _{3} 
- \eta _{7} \wedge\theta _{4} 
- \eta _{2} \wedge\theta _{11} 
+ (\eta _{5} - \eta _{1} - \eta _{4}) \wedge\theta _{12} 
\] \[\fl \quad \quad \quad
+ 2\,(\eta _{6} - \theta _{34}) \wedge\theta _{13} 
+ \theta _{33} \wedge\theta _{14} 
- \eta _{3} \wedge\theta _{22} 
- \theta _{44} \wedge\theta _{23} 
- \xi^{1} \wedge\eta _{21} 
- \xi^{2} \wedge\eta _{19} 
\] \[\fl \quad \quad \quad
- \xi^{3} \wedge\eta _{15}
- \xi^{4} \wedge\eta _{17},
\]
\[\fl
d \theta _{13}=
\eta _{12} \wedge\theta _{3} 
- \eta _{10} \wedge\theta _{4} 
+ {\textstyle \frac {2}{3}} \,(\eta _{5} - \eta _{1} -\eta _{4}) \wedge\theta _{13} 
- \eta _{2} \wedge\theta _{14} 
- \eta _{3} \wedge\theta _{23} 
- \theta _{44} \wedge\theta _{33} 
\] \[\fl \quad \quad \quad
+\eta _{6} \wedge\theta _{34} 
- \xi^{1} \wedge\eta _{17} 
- \xi^{2} \wedge\eta _{15}
- \xi^{3} \wedge\eta _{7} 
- \xi^{4} \wedge\eta _{11},
\]
\[\fl
d \theta _{14}=
\eta _{14} \wedge\theta _{3} 
- \eta _{12} \wedge\theta _{4} 
- 2\,\eta _{3} \wedge\theta _{13} 
+{\textstyle \frac {1}{3}} \,(\eta _{4} - 5 \,\eta _{1} + 2 \,\eta _{5}) \wedge\theta _{14}
+ (\eta _{6} +\theta _{34})\wedge\theta _{44}
\] \[\fl \quad \quad \quad
- \xi^{1} \wedge\eta _{18} 
- \xi^{2} \wedge\eta _{17} 
- \xi^{3} \wedge\eta _{11} 
- \xi^{4} \wedge\eta _{13},
\]
\[\fl
d \theta _{22}=
- \eta _{9} \wedge\theta _{1} 
+ \eta _{10} \wedge\theta _{2} 
+ \eta _{7} \wedge\theta _{3} 
- \eta _{8} \wedge\theta _{4} 
- 2\,\eta _{2} \wedge\theta _{12} 
+ 2\,\theta _{33} \wedge\theta _{13} 
- \xi^{1} \wedge\eta _{19} 
\] \[\fl \quad \quad \quad
+(\eta_5- 2\,\eta _{4}) \wedge\theta _{22} 
+ 2\,(\eta _{6} - 2\,\theta _{34}) \wedge\theta _{23} 
- \xi^{2} \wedge\eta _{20} 
- \xi^{3} \wedge\eta _{16}
- \xi^{4} \wedge\eta _{15},
\]
\[\fl
d \theta _{23}= 
\eta _{10} \wedge\theta _{3} 
- \eta _{9} \wedge\theta _{4} 
- 2\,\eta _{2} \wedge\theta _{13} 
+ {\textstyle \frac {1}{3}} \,(\eta _{1} - 5 \,\eta _{4}+ 2 \,\eta _{5}) \wedge\theta _{23}
+ (\eta _{6} - 3\,\theta _{34}) \wedge\theta _{33} 
\] \[\fl \quad \quad \quad
- \xi^{1} \wedge\eta _{15}
- \xi^{2} \wedge\eta _{16} 
- \xi^{3} \wedge\eta _{8} 
- \xi^{4} \wedge\eta _{7}, 
\]
\[\fl
d \theta _{33}=
{\textstyle \frac {1}{3}} \,(\eta _{5} +2 \,\eta _{1}-4\,\eta _{4}) \wedge\theta _{33} 
- 2\,\eta _{2} \wedge\theta _{34} 
- \xi^{1} \wedge\eta _{7} 
- \xi^{2} \wedge\eta _{8}
- \xi^{3} \wedge\eta _{9} 
- \xi^{4} \wedge\eta _{10},
\]
\[\fl
d \theta _{34}= 
- \eta _{3} \wedge\theta _{33} 
+ {\textstyle \frac {1}{3}} \,(\eta _{5} - \eta _{1} - \eta _{4}) \wedge\theta _{34} 
- \eta _{2} \wedge\theta _{44}
- \xi^{1} \wedge\eta _{11} 
- \xi^{2} \wedge\eta _{7} 
- \xi^{3} \wedge\eta _{10} 
\] \[\fl \quad \quad \quad
- \xi^{4} \wedge\eta _{12},
\]
\[\fl
d \theta _{44}= 
- 2\,\eta _{3} \wedge\theta _{34} 
+ {\textstyle \frac {1}{3}} \,(\eta _{5}- 4\,\eta _{1} + 2 \,\eta _{4}) \wedge\theta _{44} 
- \xi^{1} \wedge\eta _{13}
- \xi^{2} \wedge\eta _{11} 
- \xi^{3} \wedge\eta _{12} 
- \xi^{4} \wedge\eta _{14},
\]
\[\fl
d \eta _{1}= 
\eta _{2} \wedge\eta _{3}
- \eta _{12} \wedge\xi^{1} 
- \eta _{10} \wedge\xi^{2},
\]
\[\fl
d \eta _{2}= 
(\eta _{1}- \eta _{4}) \wedge\eta _{2} 
- \eta _{10} \wedge\xi^{1} 
- \eta _{9} \wedge\xi^{2},
\]
\[\fl
d \eta _{3}=
(\eta _{4} -\eta _{1}) \wedge\eta _{3}
+\eta _{14} \wedge\xi^{1} + \eta _{12} \wedge\xi^{2},
\]
\[\fl
d \eta _{4}=
- \eta _{2} \wedge\eta _{3}
+\eta _{12} \wedge\xi^{1} 
+ \eta _{10} \wedge\xi^{2},
\]
\[\fl
d \eta _{5}=0,
\]
\[\fl
d \eta _{6}=
{\textstyle \frac {1}{3}} \,(\eta _{5} - \eta _{1}-\eta_4) \wedge\eta _{6} 
-\eta_3\wedge\theta _{33} 
- \eta _{2} \wedge\theta _{44} 
+\eta _{11} \wedge\xi^{1} 
+ \eta _{7} \wedge\xi^{2} 
+ \eta _{10} \wedge\xi^{3} 
\] \[\fl \quad \quad \quad
+ \eta _{12} \wedge\xi^{4} ,
\]
\[\fl
d \eta _{7}=
{\textstyle \frac {1}{3}} \,(\eta _{5} - \eta _{1}- \eta _{4}) \wedge\eta _{7} 
- 2\,\eta _{2} \wedge\eta _{11} 
- \eta _{3} \wedge\eta _{8} 
+ \eta _{6} \wedge\eta _{10} 
- 2\,\eta _{12} \wedge\theta _{33} 
+ 2\,\eta_{10}\wedge\theta _{34}
\] \[\fl \quad \quad \quad
+\eta_9 \wedge\theta _{44} 
+\eta _{23} \wedge\xi^{1}  
+ \eta _{24} \wedge\xi^{2}
+ \eta _{25} \wedge\xi^{3} 
+ \eta _{26} \wedge\xi^{4}, 
\]
\[\fl
d \eta _{8}=
{\textstyle \frac {1}{3}} \,(\eta _{5} +2\,\eta _{1} -7 \,\eta _{4}) \wedge\eta _{8} 
- 3\,\eta _{2} \wedge\eta _{7} 
+ \eta _{6} \wedge\eta _{9} 
- 3\,\eta_{10} \wedge\theta _{33}
+ 4\,\eta _{9} \wedge\theta _{34} 
+ \eta _{24} \wedge\xi^{1}
\] \[\fl \quad \quad \quad
+\eta _{27} \wedge\xi^{2} 
+ \eta _{28} \wedge\xi^{3} 
+ \eta _{25} \wedge\xi^{4},
\]
\[\fl
d \eta _{9}=
(\eta _{1} - 2\,\eta _{4}) \wedge\eta _{9} 
- 3\,\eta _{2} \wedge\eta _{10}
+ \eta _{25} \wedge\xi^{1} 
+\eta _{28} \wedge\xi^{2},
\]
\[\fl
d \eta _{10}=
- 2\,\eta _{2} \wedge\eta _{12} 
- \eta _{3} \wedge\eta _{9} 
- \eta _{4} \wedge\eta _{10} 
+\eta _{26} \wedge\xi^{1} 
+ \eta _{25} \wedge\xi^{2},
\]
\[\fl
d \eta _{11}=
{\textstyle \frac {1}{3}} \,(\eta _{5} - \eta _{1} - \eta _{4}) \wedge\eta _{11} 
- \eta _{2} \wedge\eta _{13} 
- 2\,\eta _{3} \wedge\eta _{7} 
+ \eta _{6} \wedge\eta _{12} 
+\eta _{29} \wedge\xi^{1} 
+ \eta _{23} \wedge\xi^{2} 
\] \[\fl \quad \quad \quad
+ \eta _{26} \wedge\xi^{3}
+ \eta _{30} \wedge\xi^{4} 
+ \eta_{14}\wedge\theta _{33}
+ 2\,\eta _{10} \wedge\theta _{44},
\]
\[\fl
d \eta _{12}= 
- \eta _{1} \wedge\eta _{12} 
- \eta _{2} \wedge\eta _{14}
- 2\,\eta _{3} \wedge\eta _{10}
 + \eta _{30} \wedge\xi^{1} 
+ \eta _{26} \wedge\xi^{2},
\]
\[\fl
d \eta _{13}=
{\textstyle \frac {1}{3}} \,(\eta _{5} - 7\,\eta _{1} + 2 \,\eta _{4}) \wedge \eta _{13} 
- 3\,\eta _{3} \wedge\eta _{11} 
+ (\eta _{6}+2\,\theta_{34}) \wedge\eta _{14} 
+ 3\,\eta _{12} \wedge\theta _{44} 
+\eta _{31} \wedge\xi^{1} 
\] \[\fl \quad \quad \quad
+ \eta _{29} \wedge\xi^{2} 
+ \eta _{30} \wedge\xi^{3}
+ \eta _{32} \wedge\xi^{4}, 
\]
\[\fl
d \eta _{14}=
(\eta _{4} - 2\,\eta _{1} )\wedge\eta _{14} 
- 3\,\eta _{3} \wedge\eta _{12} 
+\eta _{32} \wedge\xi^{1} 
+ \eta _{30} \wedge\xi^{2}, 
\]
\[\fl
d \eta _{15}=
{\textstyle \frac {1}{3}} (2\,\eta _{5} - 5\,\eta _{4} - 2 \,\eta _{1}) \wedge\eta _{15} 
- 2\,\eta _{2} \wedge\eta _{17} 
- \eta _{3} \wedge\eta _{16} 
+ 2\,\eta _{6} \wedge\eta _{7} 
+ \eta _{26} \wedge\theta _{3}
- \eta _{25} \wedge\theta _{4} 
\] \[\fl \quad \quad \quad
+ \eta _{10} \wedge\theta _{13} 
+ \eta _{9} \wedge\theta _{14} 
- 2\,\eta _{12} \wedge\theta _{23} 
- 2\,\eta _{11} \wedge\theta _{33}
+ 3\,\eta _{7} \wedge\theta _{34} 
+ \eta _{8} \wedge\theta _{44} 
+ \eta _{33} \wedge\xi^{1}
\] \[\fl \quad \quad \quad
+ \eta _{34} \wedge\xi^{2} 
+\eta _{24} \wedge\xi^{3} 
+ \eta _{23} \wedge\xi^{4},
\]
\[\fl
d \eta _{16}=
{\textstyle \frac {1}{3}} \,(\eta _{1} - 8\,\eta _{4}+ 2 \,\eta _{5}) \wedge\eta _{16} 
- 3\,\eta _{2} \wedge\eta _{15} 
+ 2\,\eta _{6} \wedge\eta _{8} 
+ \eta _{25} \wedge\theta _{3} 
- \eta _{28} \wedge\theta _{4} 
+ 3\,\eta _{9} \wedge\theta _{13} 
\] \[\fl \quad \quad \quad
- 3\,\eta _{10} \wedge\theta _{23} 
+ 3\,\eta_7\wedge\theta _{33} 
+ 5\,\eta _{8} \wedge\theta _{34} 
+ \eta _{34} \wedge\xi^{1} 
+\eta _{35} \wedge\xi^{2} 
+ \eta _{27} \wedge\xi^{3}
+ \eta _{24} \wedge\xi^{4},
\]
\[\fl
d \eta _{17}=
{\textstyle \frac {1}{3}} \,(2\,\eta _{5} - 5\,\eta _{1} - 2 \,\eta _{4}) \wedge\eta _{17} 
- \eta _{2} \wedge\eta _{18} 
- 2\,\eta _{3} \wedge\eta _{15} 
+ 2\,\eta _{6} \wedge\eta _{11} 
+ \eta _{30} \wedge\theta _{3} 
\] \[\fl \quad \quad \quad
- \eta _{26} \wedge\theta _{4} 
- \eta _{12} \wedge\theta _{13} 
+ 2\,\eta _{10} \wedge\theta _{14} 
- \eta _{14} \wedge\theta _{23} 
-\eta_{13}\wedge\theta _{33}
+ \eta _{11} \wedge\theta _{34} 
+ 2\,\eta_7 \wedge\theta _{44} 
\] \[\fl \quad \quad \quad
+ \eta _{36} \wedge\xi^{1} 
+ \eta _{33} \wedge\xi^{2} 
+ \eta _{23} \wedge\xi^{3} 
+ \eta _{29} \wedge\xi^{4}, 
\]
\[\fl
d \eta _{18}=
{\textstyle \frac {1}{3}} \,(\eta _{4} - 8 \,\eta _{1}+ 2 \,\eta _{5}) \wedge\eta _{18} 
- 3\,\eta _{3} \wedge\eta _{17} 
+ 2\,\eta _{6} \wedge\eta _{13} 
+ \eta _{32} \wedge\theta _{3} 
- \eta _{30} \wedge\theta _{4} 
\] \[\fl \quad \quad \quad
- 3\,\eta _{14} \wedge\theta _{13} 
+ 3\,\eta _{12} \wedge\theta _{14} 
-\eta_{13} \wedge\theta _{34}
+ 3\,\eta_{11} \wedge\theta _{44} 
+ \eta _{37} \wedge\xi^{1} 
+ \eta _{36} \wedge\xi^{2} 
+ \eta _{29} \wedge\xi^{3} 
\] \[\fl \quad \quad \quad
+\eta _{31} \wedge\xi^{4}, 
\]
\[\fl
d \eta _{19}= 
(\eta _{5} - \eta _{1}- 2\,\eta _{4}) \wedge\eta _{19} 
- 2\,\eta _{2} \wedge\eta _{21} 
- \eta _{3} \wedge\eta _{20} 
+ 3\,\eta _{6} \wedge\eta _{15} 
- \eta _{25} \wedge\theta _{1} 
+ \eta _{26} \wedge\theta _{2} 
\] \[\fl \quad \quad \quad
+ \eta _{23} \wedge\theta _{3} 
- \eta _{24} \wedge\theta _{4}
+ \eta _{9} \wedge\theta _{11} 
+ \eta _{10} \wedge\theta _{12} 
+ \eta _{7} \wedge\theta _{13} 
+ \eta _{8} \wedge\theta _{14} 
- 2\,\eta _{12} \wedge\theta _{22} 
\] \[\fl \quad \quad \quad
- 2\,\eta _{11} \wedge\theta _{23} 
- 2\,\eta_{17} \wedge\theta _{33}
+ 4\,\eta_{15} \wedge\theta _{34}
+ \eta_{16} \wedge\theta _{44}
+\eta _{38} \wedge\xi^{1} 
+ \eta _{39} \wedge\xi^{2} 
+ \eta _{34} \wedge\xi^{3} 
\] \[\fl \quad \quad \quad
+ \eta _{33} \wedge\xi^{4}, 
\]
\[\fl
d \eta _{20}=
(\eta _{5} - 3\,\eta _{4}) \wedge\eta _{20} 
- 3\,(\eta _{2} \wedge\eta _{19} - \eta _{6} \wedge\eta _{16}) 
- \eta _{28} \wedge\theta _{1} 
+ \eta _{25} \wedge\theta _{2} 
+ \eta _{24} \wedge\theta _{3} 
\] \[\fl \quad \quad \quad
- \eta _{27} \wedge\theta _{4}
+ 3\,(\eta _{9} \wedge\theta _{12} 
+ \eta _{8} \wedge\theta _{13} 
- \eta _{10} \wedge\theta _{22} 
- \eta _{7} \wedge\theta _{23} 
- \eta_{15} \wedge\theta _{33} 
+ 2\,\eta_{16} \wedge\theta _{34})
\] \[\fl \quad \quad \quad
+ \eta _{39} \wedge\xi^{1} 
+ \eta _{40} \wedge\xi^{2} 
+ \eta _{35} \wedge\xi^{3}
+ \eta _{34} \wedge\xi^{4},
\]
\[\fl
d \eta _{21}=
(\eta _{5} - 2\,\eta _{1}- \eta _{4}) \wedge\eta _{21} 
- \eta _{2} \wedge\eta _{22} 
- 2\,\eta _{3} \wedge\eta _{19} 
+ (3\,\eta _{6} - 2\,\theta _{34}) \wedge\eta _{17} 
- \eta _{26} \wedge\theta _{1} 
\] \[\fl \quad \quad \quad
+ \eta _{30} \wedge\theta _{2} 
+ \eta _{29} \wedge\theta _{3} 
- \eta _{23} \wedge\theta _{4} 
+ 2\,\eta _{10} \wedge\theta _{11} 
- \eta _{12} \wedge\theta _{12} 
- \eta _{11} \wedge\theta _{13} 
+ 2\,\eta _{7} \wedge\theta _{14}  
\] \[\fl \quad \quad \quad
- \eta _{14} \wedge\theta _{22} 
- \eta _{13} \wedge\theta _{23} 
- \eta_{18} \wedge\theta _{33}
+ 2\,\eta_{15} \wedge\theta _{44}
+ \eta _{41} \wedge\xi^{1} 
+ \eta _{38} \wedge\xi^{2} 
+ \eta _{33} \wedge\xi^{3} 
\] \[\fl \quad \quad \quad
+ \eta _{36} \wedge\xi^{4},
\]
\[\fl
d \eta _{22}=
(\eta _{5} - 3\,\eta _{1}) \wedge\eta _{22} 
- 3\,\eta _{3} \wedge\eta _{21} 
+ 3\,\eta _{6} \wedge\eta _{18} 
- \eta _{30} \wedge\theta _{1}
+ \eta _{32} \wedge\theta _{2} 
+ \eta _{31} \wedge\theta _{3}
\] \[\fl \quad \quad \quad
- \eta _{29} \wedge\theta _{4} 
+ 3\,(\eta _{12} \wedge\theta _{11} 
- \eta _{14} \wedge\theta _{12} 
- \eta _{13} \wedge\theta _{13} 
+ \eta _{11} \wedge\theta _{14} 
+ \eta_{17} \wedge\theta _{44}
)
\] 
\begin{equation}
\fl \quad \quad \quad
+\eta _{42} \wedge\xi^{1} 
+ \eta _{41} \wedge\xi^{2} 
+ \eta _{36} \wedge\xi^{3} 
+ \eta _{37} \wedge\xi^{4}.
\label{SEs}
\end{equation}
For these equations, the non-zero reduced Cartan's characters are $s^{\prime}_1 = 16$ and $s^{\prime}_2 = 4$,  the degree of indeterminancy is $r^{(2)} = 24$, therefore Eqs. (\ref{SEs}) are  involutive, and dif\-fe\-o\-mor\-phisms from the symmetry pseudo-group depend on 4 arbitrary functions of two variables.

 In the next calculations  we use the following Maurer--Cartan forms only:  
\[\fl
\theta_0 = b_{3}^3\,b_0\,\vartheta_0,
\]
\[\fl
\theta_1 = 
b_{3}^3\,\left(
b_{22}\,\vartheta_1
- b_{21}\,\vartheta_2
+(b_{22}\,u_{zz}-b_{21}\,(u_{yz}+b_{4}))\,\vartheta_3
+(b_{21}\,u_{yy}-b_{22}\,(u_{yz}-b_{4}))\,\vartheta_4
\right),
\]
\[\fl
\theta_2 = 
b_{3}^3\,\left(
-b_{12}\,\vartheta_1
+ b_{11}\,\vartheta_2
+(b_{11}\,(u_{yz}+b_{4})-b_{12}\,u_{zz})\,\vartheta_3
+(b_{12}\,(u_{yz}-b_{4})-b_{11}\,u_{yy})\,\vartheta_4
\right),
\]
\[\fl
\theta_3 = 
b_{3}^2\,\left(
-b_{11}\,\vartheta_3+b_{12}\,\vartheta_4
\right),
\]
\[\fl
\theta_4 = 
b_{3}\,\left(
b_{21}\,\vartheta_3-b_{22}\,\vartheta_4
\right),
\]
\[\fl
\xi^1 = b_{11}\,dt+b_{12}\,dx,
\]
\[\fl
\xi^2 = b_{21}\,dt+b_{22}\,dx,
\]
\[\fl
\xi^3 = b_{3}\,\left(
(b_{22}\,u_{zz} - b_{21}\,(u_{yz}-b_{4}))\,dt
+
(b_{22}\,(u_{yz}-b_{4})-b_{21}\,u_{yy})\,dx
-b_{22}\,dy-b_{21}\,dz
\right),
\]
\[\fl
\xi^4 = b_{3}\,\left(
(b_{11}\,(u_{yz}-b_{4}) - b_{12}\,u_{zz})\,dt
+
(b_{12}\,(u_{yz}+b_{4})-b_{11}\,u_{yy})\,dx
-b_{12}\,dy-b_{11}\,dz
\right),
\]
\[\fl
\theta_{33} = \frac{b_{3}}{b_0}\,\left(b_{11}^2\,\bar{\vartheta}_{33}-2\,b_{11}b_{12}\,\bar{\vartheta}_{34}+b_{12}^2\,\bar{\vartheta}_{44}\right),
\]
\[\fl
\theta_{34} = -\frac{b_{3}}{b_0}\,\left(b_{11}\,b_{21}\,\bar{\vartheta}_{33}
-(b_{11}b_{22}+b_{12}b_{21})\,\bar{\vartheta}_{34}+b_{12}b_{22}\,\bar{\vartheta}_{44}\right),
\]
\[\fl
\theta_{44} =\frac{b_{3}}{b_0}\,\left(b_{21}^2\,\bar{\vartheta}_{33}
-2\,b_{21}b_{22}\,\bar{\vartheta}_{34}+b_{22}^2\,\bar{\vartheta}_{44}\right),
\]
\[\fl
\eta_1 = \frac{1}{b_0} \,(b_{22}\,db_{11}-b_{21}\,db_{12})
-\frac{1}{b_0^2}\,\left(
(b_{11}b_{21}^2\,u_{yyy}-b_{21}\,(2\,b_{11}b_{22}+b_{12}b_{21})\,u_{yyz}
-b_{12}b_{22}^2\,u_{zzz}
\right.
\] \[\fl \quad \quad \quad
\left.
+b_{22}\,(b_{11}b_{22}+2\,b_{12}b_{21})\,u_{yzz}
)\,\xi^1
-(b_{11}^2b_{21}\,u_{yyy}-b_{11}\,(b_{11}b_{22}+2\,b_{12}b_{21})\,u_{yyz}
\right.
\] \[\fl \quad \quad \quad
\left.
+b_{12}\,(2\,b_{11}b_{22}+b_{12}b_{21})\,u_{yzz}
-b_{12}^2b_{22}\,u_{zzz})\,\xi^2
\right),
\]
\[\fl
\eta_2 = \frac{1}{b_0} \,(b_{11}\,db_{12}-b_{12}\,db_{11})
+\frac{1}{b_0^2}\,\left(
(b_{11}^2b_{21}\,u_{yyy}-b_{11}\,(b_{11}b_{22}+2\,b_{12}b_{21})\,u_{yyz}
-b_{12}^2b_{22}\,u_{zzz}
\right.
\] \[\fl \quad \quad \quad
\left.
+b_{12}\,(2\,b_{11}b_{22}+b_{12}b_{21})\,u_{yzz})\,\xi^1-(b_{11}^3\,u_{yyy}-3\,b_{11}^2b_{12}\,u_{yyz}
+3\,b_{11}b_{12}^2\,u_{yzz}
\right.
\] \[\fl \quad \quad \quad
\left.
-b_{12}^3\,u_{zzz})\,\xi^2
\right),
\]
\[\fl
\eta_3 = \frac{1}{b_0}\,\left(b_{22}db_{21}-b_{22}^2db_{11}-b_{21}db_{22}+b_{21}b_{22}	db_{12}\right)
+\frac{b_{22}}{b_{12}}\,\eta_1
\] \[\fl \quad \quad \quad
+\frac{1}{b_{12}b_0}\,\left(
(b_{21}^2\,u_{yyy}-2\,b_{21}b_{22}\,u_{yyz}+b_{22}^2\,u_{yzz})\,\xi^1
\right.
\] \[\fl \quad \quad \quad
\left.
-(b_{11}b_{21}\,u_{yyy}-(b_{11}b_{22}+b_{12}b_{21})\,u_{yyz}+b_{12}b_{22}\,u_{yzz})\,\xi^2
\right),
\]
\[\fl
\eta_4 = \frac{1}{b_0}\,\left(
b_{12}b_{21}\,db_{11}-b_{21}\,db_{12}-b_{12}\,db_{21}+b_{11}\,db_{22}
\right)
+\frac{b_{21}}{b_{11}}\,\eta_2
\] \[\fl \quad \quad \quad
+\frac{1}{b_{11}b_0}\,\left(
((b_{11}b_{22}+b_{12}b_{21})\,u_{yzz}-b_{11}b_{21}\,u_{yyz}-b_{12}b_{22}\,u_{zzz})\,\xi^1
\right.
\] \[\fl \quad \quad \quad
\left.
+(b_{11}^2\,u_{yyz}-2\,b_{11}b_{12}\,u_{yzz}+b_{12}^2u_{zzz})\,\xi^2
\right),
\]
\[\fl
\eta_5 = 3\,\left(
\frac{db_{3}}{b_{3}}
+\frac{1}{b_{11}b_0}\,(b_{12}b_{21}db_{11}-b_{11}b_{21}db_{12}-b_{11}b_{12}db_{21}+b_{11}^2db_{22})
+\frac{b_{21}}{b_{11}}\,\eta_2
\right)+\eta_1
\] \[\fl \quad \quad \quad
-2\,\eta_4
+\frac{3}{b_{11}b_0}\,\left(
((b_{11}b_{22}+b_{12}b_{21})\,u_{yzz}-b_{11}b_{21}\,u_{yyz}-b_{12}b_{22}\,u_{zzz})\,\xi^1
\right.
\] \[\fl \quad \quad \quad
\left.
+(b_{11}^2\,u_{yyz}-2\,b_{11}b_{12}\,u_{yzz}+b_{12}^2\,u_{zzz})\,\xi^2
\right),
\]
\begin{equation}
\fl
\eta_6 = -b_{3}\,db_{4}+\theta_{34}.
\label{mcfs}
\end{equation}
In these forms, 
$\bar{\vartheta}_{ij} = \iota^{*}_0\,\vartheta_{ij}$  and
$b_{11}$, $b_{12}$, $b_{21}$, $b_{22}$, $b_{3}$, $b_{4}$ are arbitrary parameters such that  
$b_0= b_{11}b_{22}-b_{12}b_{21}\not =0$ and  $b_{11}\,b_{3}\not =0$.

\section{Contact integrable extensions}

To apply \'Eli Cartan's structure theory of Lie pseudo-groups to the problem of finding zero-curvature representations we use the notion of integrable extension. It was introduced in \cite{BryantGriffiths} for the case of {\sc pde}s with two independent variables and finite-di\-men\-si\-o\-nal coverings. The generalization of the definition to the case of infinite-di\-men\-si\-o\-nal coverings of {\sc pde}s with more than two independent variables is proposed in \cite{Morozov2009a}.  In contrast to \cite{WE,BryantGriffiths}, the starting point of our definition is the set of Maurer--Cartan forms of the symmetry pseudo-group of a given {\sc pde}, and all the constructions are carried out in terms of invariants of the pseudo-group. Therefore, the effectiveness of our method increases when it is applied to equations with large symmetry pseudo-groups.

Let $\mathfrak{G}$  be a Lie pseudo-group on a manifold $M$. Let  $\omega^1$, ... , $\omega^m$, $m=\mathrm{dim}\,M$,  be its Maurer--Cartan forms with the structure equations 
\begin{equation}
d \omega^i = A_{\gamma j}^i \,\pi^\gamma \wedge \omega^j + B_{jk}^i\,\omega^j \wedge \omega^k,
\label{SEs_of_a_Lie_pseudogroup} 
\end{equation}
where $\gamma \in \{1,\dots,\Gamma\}$ for some  $\Gamma \ge 0$. 
The coefficients $A_{\gamma j}^i$, $B_{jk}^i = - B_{kj}^i$ in these equations depend on the invariants $U^\kappa$, 
$\kappa \in \{1, \dots, \Lambda \}$, $\Lambda \ge 0$. The  differentials of the invariants satisfy equations
\begin{equation}
d U^\lambda = C_j^\lambda\,\omega^j,
\label{dUs_of_a_Lie_pseudogroup}
\end{equation} 
where $C_j^\lambda$ are functions of $U^\kappa$. Consider the following system of equations
\begin{eqnarray}
\fl
d\tau^q &=& 
D^q_{\rho r} \, \eta^\rho \wedge \tau^r 
+ 
E^q_{r s} \, \tau^r \wedge \tau^s 
+
F^q_{r \beta} \, \tau^r \wedge \pi^\beta
+
G^q_{r j} \, \tau^r \wedge \omega^j
+
H^q_{\beta j} \, \pi^\beta \wedge \omega^j
\nonumber
\\
\fl
&&+
I^q_{j k} \, \omega^j \wedge \omega^k,
\label{extra_SE} 
\\
\fl
d V^\epsilon &=&  J^\epsilon_j \,\omega^j 
+  K^\epsilon_q \, \tau^q,
\label{dVs}
\end{eqnarray}
for unknown 1-forms $\tau^q$, $q \in \{1,...,Q\}$, $\eta^\rho$, $\rho \in \{1,...,R\}$,
and unknown functions $V^\epsilon$, $\epsilon \in \{1,...,S\}$ with some 
$Q, R, S \in \mathbb{N}$. 
The coefficients $D^\kappa_{\rho r}$, ..., $K^\epsilon_q$ in Eqs.  (\ref{extra_SE}), (\ref{dVs}) are supposed to be fucntions of  $U^\lambda$ and $V^\gamma$.

\vskip 3 pt

\noindent
{\sc definition} 1.\,\,The system (\ref{extra_SE}), (\ref{dVs}) is called an  {\it integrable extension} of the system (\ref{SEs_of_a_Lie_pseudogroup}), (\ref{dUs_of_a_Lie_pseudogroup}), if Eqs.  (\ref{extra_SE}), (\ref{dVs}), (\ref{SEs_of_a_Lie_pseudogroup}) è (\ref{dUs_of_a_Lie_pseudogroup})
together meet the involutivity conditions and the compatibility conditions  
\begin{equation}
d(d\tau^q) \equiv 0,
\qquad 
d(d V^\epsilon) \equiv 0.
\label{compatibility_condition}
\end{equation}

\vskip 5 pt

Eqs.  (\ref{compatibility_condition}) give an over-determined system of {\sc pde}s for the coefficients 
$D^\kappa_{\rho r}$, ..., $K^\epsilon_q$ in Eqs.  (\ref{extra_SE}), (\ref{dVs}). If this system is satisfied, the third inverse fundamental Lie's theorem in Cartan's form, \cite[\S\S 16, 22--24]{Cartan1}, \cite{Cartan4}, 
\cite[\S\S 16, 19, 20, 25,26]{Vasilieva}, \cite[\S\S 14.1--14.3]{Stormark},  ensures the existence of the forms $\tau^q$, $V^\epsilon$, the solutions to Eqs. (\ref{extra_SE}), (\ref{dVs}). In acccordance with the second inverse fundamenatl Lie's theorem, the forms $\tau^q$,  $\omega^i$  are Maurer--Cartan forms for a  Lie pseudo-group  $\mathfrak{H}$  acting on  $M \times \mathbb{R}^Q$.

\vskip 3 pt

\noindent
{\sc definition} 2.\,\, The integrable extension (\ref{extra_SE}), (\ref{dVs}) is called {\it trivial}, if there exists a change of variables on the manifold of action of the pseudo-group  $\mathfrak{H}$ such that in the new coordinates the coefficients $F^q_{r \beta}$, $G^q_{r j}$, $H^q_{\beta j}$, $I^q_{j k}$ and $J^\epsilon_j$ are identically equal to zero, while the coefficients  $D^q_{\rho r}$, 
$E^q_{r s}$ and $K^\epsilon_q$ are independent of  $U^\lambda$.  Otherwise, the integrable extension is called  {\it nontrivial}.

\vskip 5pt 

Let  $\theta^\alpha_I$ and $\xi^j$ be a set of Maurer--Cartan forms of a symmetry pseudo-group $\mathfrak{Lie}(\EuScript{E})$ of a {\sc pde}  $\EuScript{E}$ such that $\xi^i$ are horizontal forms, that is, $\xi^1 \wedge ... \wedge \xi^n \not = 0$ on each solution of $\EuScript{E}$, while   $\theta_I^\alpha$ are contact forms, that is, they are equal to 0 on each solution.

\vskip 3 pt
\noindent
{\sc definition} 3.\,\,
Nontrivial integrable extension of the structure equations for the pseudo-group $\mathfrak{Lie}(\EuScript{E})$  of the form
\begin{equation}
d \omega^q =\Pi^q_r \wedge \omega^r + \xi^j \wedge \Omega^q_j, 
\label{contact_ie}
\end{equation}
$q,r \in \{1, \dots , N\}$, $N \ge 1$, is called a 
{\it contact integrable extension}, if the follownig conditions are satisfied:

\begin{enumerate}

\item[({\it i})]
$\Omega^q_j \in \langle \theta^{\alpha}_I, \, \omega^r_i\rangle_{\tt lin}$ for some additional 1-forms  $\omega^r_i$;

\item[({\it ii})]
$\Omega^q_j \not \in \langle  \omega^r_i\rangle_{\tt lin}$ for some $q$ and $j$;

\item[({\it iii})]
$\Omega^q_j \not \in \langle  \theta^{\alpha}_I \rangle_{\tt lin}$  for some $q$ and $j$;

\item[({\it iv})] 
$\Pi^q_r \in \langle \theta^{\alpha}_I, \,\xi^j, \, \omega^r, \, \omega^r_i\rangle_{\tt lin}$.

\item[({\it v})]
The coefficients of expansions of the forms $\Omega^q_j$ with respect to  $\{\theta^{\alpha}_I, \, \omega^r_i\}$
and the forms $\Pi^q_r$ with respect ot  $\{\theta^{\alpha}_I, \,\xi^j, \, \omega^r, \, \omega^r_i\}$
depend either on the invariants of the pseudo-group $\mathfrak{Lie}(\EuScript{E})$ alone, or they depend also on a set of some additional functions   $W_\rho$, $\rho \in \{1, \dots, \Lambda\}$, 
$\Lambda \ge 1$. In the latter case, there exist functions $P^{I \rho}_\alpha$, $Q^\rho_q$, $R_q^{j\rho}$  and $S_j^\rho$ such that
\begin{equation}
dW_\rho = P^{I}_{\rho \alpha}\,\theta^\alpha_I+Q_{\rho q}\,\omega^q+R_{\rho q}^{j}\,\omega^q_j+S_{\rho j}\,\xi^j,
\label{dW_in_the_definition}
\end{equation}
and the set of equations  (\ref{dW_in_the_definition}) satisfies the compatibility conditions
\begin{equation}
d(dW_\rho)  
= d \left(
P^{I}_{\rho \alpha}\,\theta^\alpha_I+Q_{\rho q}\,\omega^q+R_{\rho q}^{j}\,\omega^q_j+S_{\rho j}\,\xi^j\right)
\equiv  0.
\label{compatibility_for_dW}
\end{equation}

\end{enumerate}

\vskip 5 pt

We apply this definition to the structure equations (\ref{SEs}). We restrict our analysis to {\sc cie}s of the form 
\begin{eqnarray}
\fl
d \omega_0 
&=& 
\left(
\sum \limits_{i=0}^4 A_i \,\theta_i 
+ \sum {}^{*} B_{ij}\,\theta_{ij}
+ \sum \limits_{s=1}^{22} C_s\,\eta_s 
+ \sum \limits_{j=1}^4 D_j\,\xi^j 
+ \sum \limits_{k=1}^2 E^k\,\omega_k  
\right) \wedge \omega_0 
\nonumber
\\
\fl
&&+
\sum \limits_{k=1}^4 \left(
\sum \limits_{i=0}^4 F_{ik}\,\theta_i 
+ \sum {}^{*} G_{ijk}\,\theta_{ij}
+ \sum \limits_{m=1}^2 H^m_k\,\omega_m
\right) \wedge \xi^k,
\label{ie_main}
\end{eqnarray}
with two additional forms $\omega_1$ and $\omega_2$ mentioned in the part (i) of Definition 3.  In (\ref{ie_main}),  $\sum {}^{*}$ means summation for all $i,j \in \mathbb{N}$ such that $1\le i \le j \le 4$, $(i,j)\not = (2,4)$. 
These equations together with Eqs. (\ref{SEs}) satisfy the requirement of involutivity. 
We assume that the coefficients of (\ref{ie_main}) are either constants or functions of additional invariants   $W_\rho$ 
mentioned in the part (v) of Definition 3. In the latter case the differentials of $W_\rho$  meet the following requirement 
\begin{equation}
\fl
d W_\rho = 
\sum \limits_{i=0}^4 I_{\rho i} \,\theta_i 
+ \sum {}^{*} J_{\rho ij}\,\theta_{ij}
+ \sum \limits_{s=22}^7 K_{\rho s}\,\eta_s 
+ \sum \limits_{j=1}^4 L_{\rho j}\,\xi^j 
+ \sum \limits_{q=0}^2 M_{\rho q}\,\omega_q.
\label{dW}
\end{equation}

Defintion 3 yields an over-determined system for the coefficients of (\ref{ie_main}) and (\ref{dW}).
 The results of analysis of this system are summarized in the following theorem.

\vskip 5 pt
\noindent
{\sc Theorem }1.\,\, 
{\it
There are no {\sc cie}s (\ref{ie_main}) with constant coefficients or {\sc cie}s (\ref{ie_main}), (\ref{dW}) with one additional invariant $W_1$.  Every {\sc cie} (\ref{ie_main}), (\ref{dW}) with two additional invariants $W_1$, $W_2$ is contact-equivalent either to 
\[\fl
d\omega_0 = \left(\omega_1+W_1\,\eta_2+{\textstyle \frac{1}{3}}\,(\eta_5+ 2\,\eta_4-\eta_1)\right) \wedge \omega_0
+(W_1\,\theta_{34}-\theta_{44}+W_2\,\omega_2) \wedge \xi^1
\] 
\begin{equation}
\fl \quad \quad \quad
+(W_1\,\theta_{33}-\theta_{34}+W_2\,\omega_1) \wedge \xi^2
+\omega_1 \wedge \xi^3 + \omega_2\wedge\xi^4,
\label{omega_of_first_cie}
\end{equation}
\[\fl
dW_1 = W_1\,\omega_1-\omega_2-W_1\,\eta_1+W_1^2\,\eta_2-\eta_3+W_1\,\eta_4
+Z_1\,(\omega_0+W_2\,\xi^2+\xi^3)
\] 
\begin{equation}
\fl \quad \quad \quad
+Z_2\,(W_2\,\xi^1+\xi^4),
\label{dW_1_of_first_cie}
\end{equation}
\begin{equation}
\fl
dW_2 = \eta_6-\theta_{34}+{\textstyle \frac{1}{3}}\,W_2\,(\eta_5-\eta_1-\eta_4)
+Z_3\,(\omega_0+W_2\,\xi^2+\xi^3) 
+Z_4\,(W_2\,\xi^1+W_1\,\xi^4) 
\label{dW_2_of_first_cie}
\end{equation}
or to
\[\fl
d\omega_0 = \left(\omega_2+W_1\,\eta_3+{\textstyle \frac{1}{3}}\,(\eta_5 +2\,\eta_1-\eta_4)\right) \wedge \omega_0
+(\theta_{34}-W_1\,\theta_{44}+W_2\,\omega_2) \wedge \xi^1
\] 
\begin{equation}
\fl \quad \quad \quad
+(\theta_{33}-W_1\,\theta_{34}+W_2\,\omega_1) \wedge \xi^2
+\omega_1 \wedge \xi^3 + \omega_2\wedge\xi^4,
\label{omega_of_second_cie}
\end{equation}
\[\fl
dW_1 = W_1\,\omega_2-\omega_1+W_1\,\eta_1-\eta_2+W_1^2\,\eta_3-W_1\,\eta_4
+Z_1\,(\omega_0+W_2\,\xi^1+\xi^4)
\] 
\begin{equation}
\fl \quad \quad \quad
+Z_2\,(W_2\,\xi^1+\xi^4),
\label{dW_1_of_second_cie}
\end{equation}
\begin{equation}
\fl
dW_2 = \eta_6-\theta_{34}+{\textstyle\frac{1}{3}}\,W_2\,(\eta_5-\eta_1-\eta_4)
+Z_3\,(\omega_0+W_2\,\xi^1+\xi^4) 
+Z_4\,(W_2\,\xi^2+W_1\,\xi^3),
\label{dW_2_of_second_cie}
\end{equation}
where $Z_1$, ... , $Z_4$ are arbitrary parameters.}

\vskip 5 pt

The forms (\ref{mcfs}) in Eqs. (\ref{omega_of_first_cie}), (\ref{dW_1_of_first_cie}), (\ref{dW_2_of_first_cie})
and Eqs. (\ref{omega_of_second_cie}), (\ref{dW_1_of_second_cie}), (\ref{dW_2_of_second_cie}) are known explicitly, therefore, in accordance with the third inverse fundamental Lie's theorem, the forms $\omega_0$ satisfying  (\ref{omega_of_first_cie}) or (\ref{omega_of_second_cie}) can by found by means of integration.  This analysis splits into two cases --- when $Z_3 = 0$ or $Z_3 \not =0$.  

\vskip 5 pt
\noindent
{\sc remark} 1. \,\,When $Z_3 = 0$ in Eq. (\ref{dW_2_of_first_cie}) or Eq. (\ref{dW_2_of_second_cie}), the functions $W_2$ appear to be in\-de\-pen\-dent of the fibre coordinates of the covering. This entails that one symmetry of Eq. (\ref{PlebanskiEquation}) is unliftable to the fibre of the covering. From results of \cite{KV89,Krasilshchik2000,IgoninKrasilshchik,IgoninKerstenKrasilshchik,Marvan2002} 
it follows  that the corresponding covering has a non-removable parameter. Thus the appearance of the non-removable parameter in the covering can be deduced from the form of the {\sc cie} directly, before  integration of its equations.

\vskip 5 pt
The results of integration of Eqs. (\ref{omega_of_first_cie}), (\ref{dW_1_of_first_cie}), (\ref{dW_2_of_first_cie}) and Eqs. (\ref{omega_of_second_cie}), (\ref{dW_1_of_second_cie}), (\ref{dW_2_of_second_cie}) are given in the following theorem.

\vskip 5 pt
\noindent
{\sc theorem} 2.\,\,
{\it 
When $Z_3=0$, every solution to Eq. (\ref{omega_of_first_cie}) up to a contact equivalence is 
\[
\fl
\omega_0 = \frac{b_0b_{3}}{b_{12}\,v_z - b_{11}\,v_y}\,
\left(dv
+(v_{zz}\,v_{y}- (u_{yz}+\lambda)\,v_z)\,dt
+((u_{yz}-\lambda)\,v_{y}-u_{yy}\,v_{z})\,dx
\right.
\]
\begin{equation}
\fl \quad \quad \quad
\left.
-v_y\,dy - v_z\,dz\right),
\label{first_WE_form}
\end{equation}
whereas for $Z_3 \not =0$ it is
\[
\fl
\omega_0 = \frac{b_0b_{3}}{b_{12}\,v_z - b_{11}\,v_y}\,
\left(dv
+(v_{zz}\,v_{y}- (u_{yz}+v)\,v_z)\,dt
+((u_{yz}-v)\,v_{y}-u_{yy}\,v_{z})\,dx
\right.
\]
\begin{equation}
\fl \quad \quad \quad
\left.
-v_y\,dy - v_z\,dz\right).
\label{second_WE_form}
\end{equation}

The solutions to Eq. (\ref{omega_of_second_cie}) can be obtained from (\ref{first_WE_form}) and 
(\ref{second_WE_form}) by the following simple change of independent variables:   $(t,x,y,z) \mapsto (x,t,z,y)$.
}

\vskip 5 pt
 When we put $\omega_0=0$, Eq. (\ref{first_WE_form}) gives the system (\ref{covering}), while Eq. 
(\ref{second_WE_form}) defines new covering 
\[
\left\{
\begin{array}{l}
v_t = (u_{yz} + v) \, v_z - u_{zz}\,v_y,
\\
v_x = u_{yy} \, v_z - (u_{yz} - v)\,v_y
\end{array}
\right.
\]
for the second heavenly equation. These equations are nonlinear w.r.t. the fibre variable $v$.

\vskip 5 pt
\noindent
{\sc remark} 2.\,\,Direct computation shows that the symmetry of Eq. (\ref{PlebanskiEquation}) with the in\-fi\-ni\-te\-si\-mal generator
$X = t\,\frac{\partial}{\partial y}+x\,\frac{\partial}{\partial z}$ is unliftable to a symmetry of Eqs. (\ref{covering}).
Since $e^{\lambda X} (u_{yy}, u_{yz}, u_{zz}) = (u_{yy}, u_{yz}, u_{zz})$ and
$e^{\lambda X} (v,v_t,v_x,v_y,v_z) = (v,v_t+\lambda\,v_y,v_x+\lambda\,v_z,v_y,v_z)$, the parameter $\lambda$ in 
Eqs. (\ref{covering}) can be obtained by the action of $e^{\lambda X}$ to the system (\ref{covering}) with $\lambda=0$.
Therefore, $\lambda$ is the non-removable parameter of the covering (\ref{covering}).

\end{document}